\documentclass{amsart}
\usepackage{amsmath}
\usepackage{graphicx}
\usepackage{amsfonts}
\usepackage{amssymb}

\newtheorem{theorem}{Theorem}
\theoremstyle{plain}

\newtheorem{corollary}{Corollary}

\newtheorem{lemma}{Lemma}

\newtheorem{proposition}{Proposition}

\numberwithin{equation}{section}
\newcommand{\leftexp}[2]{{\vphantom{#2}}^{#1}{#2}}

\begin{document}
\title{One dimensional conformal metric flow II}
\author{Yilong Ni and Meijun Zhu}
\address{Department of Mathematics\\
The University of Oklahoma\\
Norman, OK 73019\\
}

\begin{abstract}
In this paper we continue our studies of  the one dimensional
conformal metric flows, which were introduced in \cite{NZ1}. In
this part we mainly focus on evolution equations involving fourth
order derivatives.
 The global existence and exponential convergence of metrics for the $1$-$Q$ and
$4$-$Q$ flows are obtained.
\end{abstract}
\maketitle

\section{Introduction}
In \cite{NZ1} we initiated our study of one dimensional conformal
curvature problem. Our research  revealed the rich conformal
structures on $S^1$. It also has impacts in the study of affine
geometry and its application to image processing.

Recall that if $(S^1,g_s)$ is the unit circle with the induced
metric $g_s=d \theta \otimes d \theta$ from ${\bf R}^2$, for any
metric $g$ on $S^1$ (for example, this metric could be given by
reparametrizing the circle), we write $g:=d \sigma \otimes d
\sigma=v^{-\frac43}g_s$ for some positive function $v$. We then
introduce a general {\it $\alpha$- scalar curvature} of $g$ for
any positive constant $\alpha$  by
$$
\leftexp{\alpha}{R}_g = v(\alpha
(v^{\frac13})_{\theta\theta}+v^{\frac13}).
$$
The $\alpha-$scalar curvature flows for $\alpha=1$ and $4$ were
studied in \cite{NZ2}, where the exponential convergence of
metrics were obtained.

 We further  define a general {\it $\alpha$-$Q$ curvature} of
$g$ for any positive constant $\alpha$ by
$$
\leftexp{\alpha}{Q_g}=v^{\frac53}(\frac{\alpha^2}9 v_{\theta\theta\theta\theta}
+\frac{10\alpha}9v_{\theta\theta}+v).
$$
Thus $\leftexp{\alpha}{Q}_{g_s}=1$. The corresponding {\it
$\alpha$-conformal $\leftexp{\alpha}{P}_g$} operator of $g$ is
defined by
$$
\leftexp{\alpha}{P}_g f=\frac{\alpha^2}9\Delta_g^2f
+\frac{10\alpha}9\nabla_g(\leftexp{\alpha}{R}_g\nabla_g
f)+\leftexp{\alpha}{Q}_g f,
$$
where $\nabla_g=D_\sigma$, $\Delta_g=D_{\sigma \sigma}$ and
$\leftexp{\alpha}{R}_g$ is the $\alpha$-scalar curvature of $g$.

We shall suppress the superscript ``$\alpha$" if no confusion
would result. It is proved in \cite{NZ1} that $P_g$ is a conformal
covariant.
\begin{proposition}
If $g_2=\varphi^{-\frac43}g_1$, then
$Q_{g_2}=\varphi^{\frac53}P_{g_1}\varphi$ and
$P_{g_2}\psi=\varphi^{\frac53}P_{g_1}(\psi\varphi)$,
for any $\psi\in C^{4}({S^1})$.
\label{prop0-1}
\end{proposition}
The general {\it $\alpha$-$Q$ curvature} flow is introduced as
\begin{equation}
\partial_t g=\leftexp{\alpha}{Q}_g g.
\label{0-2}
\end{equation}
We will see in Section 2 that it is equivalent to the normalized
$\alpha$-$Q$ curvature flow:
\begin{equation}\label{normal}
\partial_t g=(\leftexp{\alpha}{Q}_g-\leftexp{\alpha}{\overline{Q}}_g)g,
\qquad L(0)=2\pi,
\end{equation}
where $\leftexp{\alpha}{\overline{Q}}_g=\int \leftexp{\alpha}{Q}_g
d\sigma \big/\int d\sigma$ and $L(t)=\int d\sigma$. It will be
clear that this flow is in fact the gradient flow of total
curvature  $\leftexp{\alpha}{\overline{Q}}_g$ (see Lemma
\ref{lem2} below).  We  pointed out in \cite{NZ1} that two cases
of $\alpha=1$ and $\alpha=4$ are of special interest. In this
paper we shall focus on these two cases. We will prove the global
existence of the flows and exponential convergence of metrics for
these two flows.

Recently there are some beautiful results on $Q$-curvature flow
equations, though all of them focus on higher dimensional cases
(for dimension $n \ge 4$). The global existence and convergence of
higher order flow on general compact manifolds were obtained by
Brendle under the condition of smaller total $Q-$curvature than
that of  the sphere with standard metric  \cite{BR}. The
convergence of $Q-$ curvature flow on $S^4$ with the initial
metric in the same conformal class of the standard metric was
later obtained in \cite{BR1}. The flow approach to the prescribing
$Q$-curvature  on $S^4$ is carried out by Malchiodi and Struwe
\cite{MS}. There are two main ingredients in the proof of global
existence and convergence of $Q$-curvature flow on $S^4$. One is
the sharp inequality involving higher order derivatives which
guarantees the lower bound for a certain functional (see, for
example, Branson, Chang and Yang \cite{BCY}, and Beckner
\cite{BE}); The other is the new approach to the flow equations
via integral estimates (see, for example, Chen \cite{CH}, and
Schwetlick and  Struwe \cite{ST}).

Even though our flow is on one dimensional circle, we face the
similar difficulty. For $\alpha=4$, the extremal metric was
classified in \cite{NZ1} (see, also Hang \cite{Ha}). We thus can
establish the global and convergence of 4-$Q$-curvature flow along
the line as we just described.

\begin{theorem}
There is a unique smooth solution $g(t)$, $t\in[0,\infty)$ to the
flow equation (\ref{normal}) for any given initial metric
$g_0=v^{-\frac43}(\theta,0)g_s$ on $S^1$. Moreover $g(t)$
converges exponentially to a smooth metric $g(\infty)$ and the
$4$-$Q$-curvature of $(S^1,g(\infty))$ is constant.
\label{theorem1}
\end{theorem}

The case of  $\alpha=1$ is more subtle. In \cite{NZ1}, we  proved
the existence of  sharp inequality (Theorem 3 in \cite{NZ1}, see
Remark 6 there for more comments), but were not able to classify
the extremal metrics. To prove the exponential convergence of
metrics under the 1-$Q$-curvature flow, one needs to classify all
extremal metrics and to know the precise sharp constant.  Now we
can achieve this:

\begin{theorem}\label{theorem2}
For $u(\theta)\in H^2(S^1)$ and $u>0$, if $u$ satisfies
\begin{equation}\label{on}
\int_0^{2\pi}\frac{\cos^{3}(\theta+\alpha)}{u^{5/3}(\theta)}d\theta=0
\end{equation}
for all $\alpha\in [0, 2\pi)$, then
$$
\int_0^{2\pi}(u^2_{\theta\theta}-10u^2_\theta+9u^2)d\theta
\left(\int_0^{2\pi}u^{-2/3}(\theta)d\theta\right)^3\ge 144\pi^4.
$$
More over, if $u_0$ is an extremal function, then the
1-$Q$-curvature of $u_0^{-4/3} g_s$ is a constant, and
\begin{equation}\label{form}
u_0(\theta)=c\left(\lambda^2\cos^2(\theta-\beta)
+\lambda^{-2}\sin^2(\theta-\beta)\right)^{\frac32},
\end{equation}
for some $ \lambda, c>0$ and $ \beta\in [0, 2\pi)$.
\end{theorem}

With this classification of extremal metrics, we are able to prove
the exponential convergence of metrics under the 1-$Q$-curvature
flow.

\begin{theorem}
Suppose the initial metric $g_0=v^{-\frac43}(\theta,0)g_s$ on $S^1$ satisfies
the orthogonal condition
$$
\int_0^{2\pi}\frac{\cos^3(\theta+\beta)}{v^{\frac{5}{3}}}d\theta=0, \quad
\mbox{ for any }\beta\in[0,2\pi).
$$
Then there is a unique solution $g(t)$ to the flow equation
(\ref{normal}) for $t\in[0,\infty)$. Moreover $g(t)$ converges
exponentially to a smooth metric $g(\infty)$ and the
$1$-$Q$-curvature of $(S^1,g(\infty))$ is constant.
\label{theorem3}
\end{theorem}

 The paper is organized as follows. In
Section 2, we derive some basic properties about the flow and
prove the global existence of the flow when $\alpha=1$ or $4$. The
$L^\infty$ convergences of the $\alpha-Q$ curvature for $\alpha
=1$ and $4$ are obtained  in Section 3. Using integral estimates,
we then prove the exponential convergence of metric for $4$-$Q$
flow in Section 4, and complete the proof of Theorem 1. In Section
5, we first classify all the constant $1$-$Q$ curvature metrics on
$S^1$, thus complete the proof of Theorem 2; Then using a similar
argument as in Section 4 we prove the exponential convergence of
metric for $1$-$Q$ flow, thereby  complete the proof of Theorem 3.

Throughout this paper, we use $C$, $C_1$, $C_2$, $\cdots$, to
represent some various positive constants.

\medskip
\noindent{\bf ACKNOWLEDGMENT. } The work of M. Zhu is partially
supported by the NSF grant DMS-0604169.

\section{Basic properties and global existence}
In this section, we shall derive some basic equations for
$\alpha$-scalar curvature, $\alpha$-Q curvature and conformal
factor function  under the flow. We then derive the {\it a priori}
$L^\infty$ estimates (depending on time) for metrics. From the
estimates we obtain the global existence for the flows.

 We  first show that the
flow equation (\ref{0-2}) for $g(t)$ is equivalent to a normalized
flow (\ref{normal}). In fact, if we choose
$$
\hat g(t)=\frac{4\pi^2}{L^2(0)}\exp \left(-\int_0^t\overline{Q}_g
(\tau) d\tau\right) g(t),
$$
where $L(t)=\int d \sigma(t)$ and $g(t)=d \sigma(t) \otimes d
\sigma(t)$, and  a new time variable
$$\hat t=\frac{16\pi^4}{L^4(0)}\int_0^t \exp\left(-2\int_0^\delta
 \overline{Q}_g (\tau)d \tau \right)d \delta,
$$
then equation (\ref{0-2}) can be written as
$$
\partial_{\hat t} \hat g= (\hat Q _{\hat g}
-\frac{\int \hat Q _{\hat g}d \hat \sigma}{\int d\hat  \sigma})
\hat g.
$$

From now on, we shall focus on the normalized flow (\ref{normal}).
\begin{lemma}
Along flow equation (\ref{normal}) with
$g(\sigma,t)=v^{-\frac43}g_s$, curvatures $R=R _g$ and $Q=Q _g$
satisfy
\begin{equation}
{R}_t= -\frac{\alpha}{4}\Delta Q-R(Q-\overline{Q})
\label{1-r}
\end{equation}
and
\begin{equation}
{Q}_t=
-\frac{\alpha^2}{12}\Delta^2Q-\frac{5\alpha}{6}\nabla(R\nabla Q)
-2Q(Q-\overline{Q}), \label{1-q}
\end{equation}
respectively.  Metric satisfies
\begin{equation}\label{1-3}
\partial_t (d \sigma)=\frac 12 (Q-\overline{Q}) d \sigma,
\end{equation}
and  $v$ satisfies
\begin{equation}
v_t=-\frac 34(Q -\overline{Q})v \label{1-4}.
\end{equation}
\label{lem1-1}
\end{lemma}
\begin{proof}
$$
(Q-\overline{Q})g=\partial_t g=(-\frac43)v^{-\frac73}v_tg_s
=(-\frac43)v^{-1}v_tg,
$$
that is $v_t=-\frac34(Q-\overline{Q})v$. Thus
$$
\partial_t (d\sigma)=(v^{-\frac23}d\theta)_t=-\frac23v^{-\frac53}v_td\theta
=\frac12(Q-\overline{Q})d\sigma.
$$
Using the conformal invariance of $L _g$ (see, for example,
\cite{NZ2}) and $P _g$, we have
\begin{align*}
R_t=&(vL_{g_s} (v^{\frac13}))_t=v_tL_{g_s} (v^{\frac13})
+vL_{g_s} (\frac13v^{-\frac23}v_t)\\
=&-\frac34(Q-\overline{Q}) R+\frac13L_{g} (\frac{v_t}{v})\\
=&-\frac{\alpha}{4}\Delta Q-R(Q-\overline{Q}),
\end{align*}
and
\begin{align*}
Q_t=&(v^{\frac53}P_{g_s} (v)_t=\frac53v^{\frac23}v_tP_{g_s} (v)
+v^{\frac53}P_{g_s} (v\frac{v_t}v)\\
=&-\frac54(Q-\overline{Q})Q+P_{g} (\frac{v_t}{v})\\
=&-\frac{\alpha^2}{12}\Delta^2Q-\frac{5\alpha}{6}\nabla(R\nabla_g Q)
-2Q(Q-\overline{Q}).
\end{align*}
\end{proof}

It follows from (\ref{1-3}) that
$$
\partial_t\int_{S^1} d\sigma=\int_{S^1}\frac12
(Q _g-\overline{Q} _g)d\sigma=0.
$$
Thus flow (\ref{normal}) preserves the arc length with respect to
metric $g$ (i.e. $\int_0^{2\pi} v^{-\frac23} d \theta=L(0)=2\pi$).
Moreover, along the flow, we see from the following lemma that the
total $Q$ curvature is strictly decreasing unless $Q _g$ is a
constant.

\begin{lemma}\label{lem2}
Along flow (\ref{normal}), we have
\begin{equation}
\partial_t \overline{Q} _g= -\frac 3{4\pi}
\int_{S^1}(Q _g-\overline{Q} _g)^2 d\sigma. \label{1-5}
\end{equation}
 \label{lem1-2}
\end{lemma}
\begin{proof}
\begin{align*}
\partial_t\overline{Q} _g=&\frac1{2\pi}\int_{S^1} (Q _g)_td\sigma
+\frac1{2\pi}\int_{S^1} Q _g\partial_t (d\sigma)\\
=&\frac1{2\pi}\int_{S^1}(-2Q _g)(Q _g-\overline{Q} _g)d\sigma
+\frac1{4\pi}\int_{S^1} Q _g(Q _g-\overline{Q} _g)d \sigma\\
=&-\frac3{4\pi}\int_{S^1} Q _g(Q _g-\overline{Q} _g)d \sigma
=-\frac3{4\pi}\int_{S^1}( Q _g-\overline{Q} _g)^2d\sigma\le0.
\end{align*}
\end{proof}

Note that (\ref{1-4}) can also be written as
\begin{equation}
v_t= -\frac{\alpha^2}{12} v^{\frac 83} \Delta_{g_s}^2 v-\frac{
5\alpha} 6 v^{\frac 83} \Delta_{g_s} v\label{1-4'}-\frac 34
v^{\frac {11}3}+\frac 34 \overline{Q} _g v.
\end{equation}

We are now ready to prove the global existence for $\alpha=4$
and $\alpha=1$.
\begin{proposition}\label{4b}
Suppose that $g(t)=v^{-\frac43}g_s$ satisfies the flow equation
(\ref{normal}) on $[0,T)$ for $\alpha=4$. Then there exists
$C=C(T)$, such that
$$
\frac{1}{C}<v(t)<C\quad\mbox{ and }\quad||v(t)||_{H^{2}}<C\quad\mbox{ on }[0,T).
$$
\end{proposition}
\begin{proof}
For any given $\lambda>0$ and $\beta \in [0, 2\pi)$, let
$$
\Psi_{\lambda,\beta}(\theta)=\left(
\lambda^2\cos^2\frac{\theta-\beta}2
+\lambda^{-2}\sin^2\frac{\theta-\beta}2\right)^{3/2},
$$

$$
\omega_{\lambda,\beta}(\theta)=\beta+\int_\beta^\theta
{\Psi_{\lambda,\beta}}^{-\frac23}d\theta = \left\{
\begin{array}{ll}
\beta+2\arctan(\lambda^{-2}\tan\frac{\theta-\beta}{2}), &\mbox{if} \ 0\le \theta -\beta \le  {\pi}\\
\beta+2\arctan(\lambda^{-2}\tan\frac{\theta-\beta}{2})+2\pi,
&\mbox{if} \ \pi<\theta -\beta \le  {2 \pi}
\end{array}
\right.
$$
and
$$
({\mathcal T}_{\lambda ,\beta}v)(\theta)
=v(\omega_{\lambda,\beta}(\theta))\Psi_{\lambda,\beta} (\theta).
$$
By Lemma 3 in \cite{NZ1}, we know for any given $t\in[0,T)$, there
exist $\lambda=\lambda(t)>0,\beta\in[0,2\pi)$ such that
$u(\theta,t):=({\mathcal T}_{\lambda,\beta}v)(\theta)$ satisfies
$\int_0^{2\pi}u(\theta,t)\cos\theta
d\theta=\int_0^{2\pi}u(\theta,t)\sin\theta d\theta=0$. For any
positive function $\varphi\in H^2(S^1)$ define functional
$$
F(\varphi)=\frac{16}9\left(\int_0^{2\pi}(\varphi_{\theta\theta}^2
-\frac52\varphi_\theta^2
+\frac9{16}\varphi^2)d\theta\right)
\left(\int_0^{2\pi}\varphi^{-\frac23}d\theta\right)^3.
$$
Suppose the Fourier expansion of $u(\theta,t)$ is
$$
u(\theta,t)=a_0+\sum_{k=2}^\infty a_k\cos(k\theta-\gamma_k).
$$
Then we have (noting that $\int_0^{2 \pi} u^{-2/3} d \theta
=2\pi$)
$$
F(u)=(2\pi)^4\left(a_0^2+\frac89\sum_{k=2}^\infty (k^4-\frac52k^2
+\frac9{16})a_k^2 \right),
$$
which implies
$$
F(u)\ge C\int_0^{2\pi}( u_{\theta\theta}^2+u^2)d\theta
$$
for some constant $C>0$. On the other hand, from the conformal
covariance of $F$, we have
$F(u(\theta,t))=F(v(\theta,t))=(2\pi)^4\overline{Q}(t) \le
(2\pi)^4\overline{Q}(0)$. It follows that $||u(t)||_{H^2}$ and
therefore $||u(t)||_{C^{1,a}}$ ($a\in(0,\frac12]$) is bounded on
$[0,T)$. Since $\int_0^{2\pi}u^{-\frac23}d\theta=2\pi$  we know
that there exists a constant $C_1$ not depending on $T$, such that
$\frac1{C_1}\le u\le C_1$ on $[0,2\pi]\times [0,T)$.

In order to obtain the estimates on $v(\theta, t)$, it suffices to
prove that $\lambda(t)$ is bounded on $[0,T)$. Suppose not, there
exists a sequence $t_i\to T$, such that $\lambda(t_i)\to\infty$.
Without loss of generality, we may assume that $\beta(t_i)\to 0$.
Then for any $\epsilon>0$,
$$
\lim_{i\to\infty}\left(\int_{-\epsilon}^\epsilon v(t_i)^{-\frac23}d\theta
+\int_{\pi-\epsilon}^{\pi+\epsilon} v(t_i)^{-\frac23}d\theta\right)=2\pi.
$$
On the other hand, for any $t>\tilde{t}$ we have
\begin{align*}
\left|\int_{-\epsilon}^\epsilon v(t)^{-\frac23}d\theta-
\int_{-\epsilon}^\epsilon v(\tilde{t})^{-\frac23}d\theta\right|
&=\left|\int_{\tilde{t}}^t \left(\int_{-\epsilon}^\epsilon
v(t)^{-\frac23}d\theta\right)_t dt\right|\\
&\le \int_{\tilde{t}}^t
\int_0^{2\pi}\frac{|Q-\overline{Q}|}2d\sigma dt\\
&\le C_2(t-\tilde{t})^{\frac12}\left(\int_{\tilde{t}}^t\int_0^{2\pi}
(Q-\overline{Q})^2d\sigma dt\right)^{\frac12}.
\end{align*}
For fixed $t\in[0,T)$,
$$
2\pi=\lim_{i\to\infty}\left(\int_{-\epsilon}^\epsilon
+\int_{\pi-\epsilon}^{\pi+\epsilon}\right)v^{-\frac23}(t_i)d\theta
\le \left(\int_{-\epsilon}^\epsilon+\int_{\pi-\epsilon}^{\pi+\epsilon}\right)
v^{-\frac23}(t)d\theta+C_3(T-t)^{\frac12}.
$$
Choosing $t$ closed enough to $T$ and choosing $\epsilon$ small
enough we get a contradiction.
\end{proof}

For $\alpha=1$ we have a similar proposition:
\begin{proposition}\label{1b}
Suppose that $g(t)=v^{-\frac43}g_s$ satisfies the flow equation
(\ref{normal}) on $[0,T)$ for $\alpha=1$. Then there exists
$C=C(T)$, such that
$$
\frac{1}{C}<v(t)<C\quad\mbox{ and }\quad||v(t)||_{H^{2}}<C\quad\mbox{ on }[0,T).
$$
\end{proposition}
\begin{proof}
For any $\lambda>0$, let
$$
\Gamma_\lambda(\theta)=\left(\lambda^2\cos^2\theta
+\lambda^{-2}\sin^2\theta\right)^{3/2}, \quad
\sigma_\lambda(\theta)=\int_0^\theta \Gamma_\lambda^{-2/3}d\theta,
$$
and
$$
({\bf T}_\lambda u)(\theta)
:=u(\sigma_\lambda(\theta))\Gamma_\lambda(\theta).
$$
Then again by Lemma 3 in \cite{NZ1} we know that for any given
$t\in[0,T)$, there exist $\lambda=\lambda(t)>0,\beta\in[0,2\pi)$
such that $u(\theta,t):=({\bf T}_{\lambda,\beta}v)(\theta)$
satisfies $\int_0^{2\pi}u(t)\cos2\theta
d\theta=\int_0^{2\pi}u(t)\sin2\theta d\theta=0$.  For any positive
function $\varphi\in H^2(S^1)$ define functional
$$
{\bf F}(\varphi)=\frac{1}9\left(\int_0^{2\pi}(\varphi_{\theta\theta}^2
-10\varphi_\theta^2
+9\varphi^2)d\theta\right)\left(\int_0^{2\pi}\varphi^{-\frac23}d\theta\right)^3.
$$
Suppose the Fourier expansion of $u(\theta,t)$ is
$$
u(\theta,t)=a_0+a_1\cos(\theta-\gamma_1)+\sum_{k=3}^\infty a_k\cos(k\theta
- \gamma_k).
$$
Then we have
$$
{\bf F}(u)=(2\pi)^4\left(a_0^2+\frac1{18}\sum_{k=4}^\infty (k^4-10k^2+9)a_k^2
\right).
$$
Integrating the nonnegative function $u(\theta)(1\pm\cos(\theta-\gamma_1))$ we
have
$$
0\le \int_0^{2\pi}u(\theta)(1\pm\cos(\theta-\gamma_1))d\theta=2\pi a_0\pm
\pi a_1,
$$
which implies $|a_1|\le 2a_0$. Similarly, integrating the
nonnegative function $u(\theta)(1\pm\cos(3\theta-\gamma_1))$ we
have$|a_3|\le 2a_0$. Hence
$$
{\bf F}(u)\ge C ||u||_{H^2}^2
$$
for some constant $C>0$. Also from the conformal covariance of
$\bf F$, we have ${\bf F}(u(\theta,t))={\bf
F}(v(\theta,t))=(2\pi)^4\overline{Q}(t) \le
(2\pi)^4\overline{Q}(0)$. It follows that $||u(t)||_{H^2}$ is
bounded on $[0,T)$. The rest of the proof will be similar to the
proof of Proposition \ref{4b}.
\end{proof}

Using a similar argument in \cite{BR} we obtain the estimates on
higher derivative of $v(\theta, t)$ as follows.
\begin{proposition}
Suppose that $g(t)=v^{-\frac43}g_s$ satisfies the flow equation
(\ref{normal}) on $[0,T)$ for $\alpha=1\mbox{ or }4$. Then
$||v(t)||_{H^{2k}(S^1)}$ is bounded on $[0,T)$ for any $k\in\bf
N$. \label{prop4}
\end{proposition}
\begin{proof}
$$
\partial_t \int_0^{2\pi}(v^{(2k)})^2d\theta
=2\int_0^{2\pi}v^{(2k)}(v_t)^{(2k)}d\theta
=-\frac32\int_0^{2\pi}v^{(2k)}((Q -\overline{Q} )v)^{(2k)}d\theta.
$$
Since $Q=\frac19v^{\frac53}(\alpha^2v^{(4)}+10\alpha v^{(2)}+9v)$,
we obtain from Proposition \ref{4b} and \ref{1b} that
\begin{align*}
\partial_t\int_0^{2\pi}(v^{(2k)})^2d\theta\le&
-\frac{\alpha^2}6\int_0^{2\pi}(v^{(2k+2)})^2v^{\frac83}d\theta
+C_1\sum_{k_1,\cdots,k_m}\int_0^{2\pi}\prod_{i=1}^m|v^{(k_i)}|
d\theta\\
\le&-C_2\int_0^{2\pi}(v^{(2k+2)})^2v^{\frac83}d\theta
+C_1\sum_{k_1,\cdots,k_m}\int_0^{2\pi}\prod_{i=1}^m|v^{(k_i)}|
d\theta,
\end{align*}
where $\Sigma$ is taken over all $m-$tuples $k_1,\cdots,k_m$ with $m\ge3$,
which satisfy $1\le k_i\le 2k+1$ and $k_1+\cdots+k_m\le 4k+4$.

For each $m-$tuple $k_1,\cdots,k_m$, let $r_i=\max\{0,\frac{k_i-\frac1{m}
-\frac32}{2k}\}$. Then we have $r:=\theta_1+\cdots+\theta_m<2$ and
$||v^{(k_i)}||_{L^{m}}\le C||v||_{H^{k_i-\frac1{m}+\frac12}}
\le C||v||_{H^2}^{1-r_i}||v||_{H^{2k+2}}^{r_i}$. It follows that
\begin{align*}
\partial_t\int_0^{2\pi}(v^{(2k)})^2d\theta\le&
-C_1\int_0^{2\pi}(v^{(2k+2)})^2d\theta
+C_2\sum_{k_1,\cdots,k_m}\int_0^{2\pi}\prod_{i=1}^m|v^{(k_i)}| d\theta\\
\le& -C_1\int_0^{2\pi}(v^{(2k+2)})^2d\theta+C_3\sum_{k_1,\cdots,k_m}
\prod_{i=1}^m||v||_{H^{k_i-\frac1{m}+\frac12}}\\
\le& -C_1\int_0^{2\pi}(v^{(2k+2)})^2d\theta+C_4\sum_{k_1,\cdots,k_m}
||v||_{H^2}^{m-r}||v||_{H^{2k+2}}^{r}\\
\le& -C_1\int_0^{2\pi}(v^{(2k+2)})^2d\theta+C_5||v||_{H^{2k+2}}
^{r}\le C_6.
\end{align*}
Hence $\int_0^{2\pi}(v^{(2k)})^2d\theta$ is bounded on $[0,T)$.
\end{proof}

For $\alpha =4$ we know that $\leftexp{\alpha}{P}_g$ is positive.
From the above proposition, we immediately get that the 4-$Q$-flow
exists on $[0,\infty)$. To show the global existence of the flow
for $\alpha=1$, we need another lemma.
\begin{lemma}
Suppose that $g(t)=v^{-\frac43}g_s$ satisfies the flow equation
(\ref{normal}) on $[0,T)$ for $\alpha=1$. If
 $\int_0^{2\pi}\cos^3 (\theta+\alpha) \cdot v^{-5/3}(\theta,0)d \theta= 0$ for all $\alpha \in [0, 2\pi)$,  then for all $t>0$,
 \begin{equation}
\int_0^{2\pi}\cos^3 (\theta+\alpha) \cdot v^{-5/3}(\theta,t)d
\theta=0\label{1-13}
 \end{equation}
  for all $\alpha \in[0, 2\pi).$
\label{lemor}
\end{lemma}

\begin{proof}
From (\ref{1-4}) and the definition of 1-$Q$ curvature we have
\begin{align*}
\partial_t \int_0^{2\pi} &\cos^3 (\theta +\alpha)\cdot v^{-5/3}(\theta,t)d
\theta\\
&=-\frac 53 \int_0^{2\pi}\cos^3 (\theta +\alpha)\cdot
v^{-8/3}(\theta,t)
v_t(\theta,t)d \theta\\
&=\frac 54 \int_0^{2\pi} \leftexp{1}Q \cos^3(\theta +\alpha)\cdot
v^{-5/3}(\theta,t)d \theta-\frac{5 \overline{\leftexp{1}Q}}4
\int_0^{2\pi}\cos^3 (\theta+\alpha) \cdot
v^{-5/3}(\theta,t)d \theta\\
&=-\frac{5 \overline{\leftexp{1}Q}}4 \int_0^{2\pi}\cos^3
(\theta+\alpha) \cdot v^{-5/3}(\theta,t)d \theta.
\end{align*}
Thus
$$
 \int_0^{2\pi}\cos^3 (\theta+\alpha) \cdot v^{-5/3}(\theta,t)d
\theta=C e^{-\int_0^t \overline{\leftexp{1}Q}(\tau) d \tau}.
$$
Since  $\int_0^{2\pi}\cos^3 (\theta+\alpha) \cdot v^{-5/3}
(\theta,0)d \theta=0$, we have $C=0$, thus $ \int_0^{2\pi}\cos^3
(\theta+\alpha) \cdot v^{-5/3}(\theta,t)d \theta=0$.
\end{proof}

For $\alpha =1$ we know  $\leftexp{\alpha}{P}_g$ is positive on
$$\{u\in H^4(S^1) \ : \ u>0, \ \int_0^{2\pi}\cos^3 (\theta+\alpha) \cdot u^{-5/3}(\theta,t)d \theta=0, \ \ \forall \alpha\in [0, 2\pi)\},$$
see, for example, Theorem 3 in \cite{NZ1}.  The global existence
of 1-$Q$-curvature flow then follows from Proposition \ref{prop4}
and Lemma \ref{lemor}.

\section{$L^\infty$ Convergence of $^\alpha Q$ along $\alpha$-$Q$-flow}
In this section, we shall follow \cite{NZ2} closely to derive the
$L^\infty$ norm convergence for curvatures. Throughout the rest of
the paper, we will only consider $\alpha=1$ or $\alpha=4$; In the
case of $\alpha=1$, we always assume the initial metric satisfies
the orthogonal condition (\ref{1-13}) (thus always satisfies
(\ref{1-13}) along the flow by Lemma \ref{lemor}); we also denote
$L^p=L^p(d\sigma)$.

For $p\ge 2$, we define
$$
G_p(t):=\int_0^{2\pi}|Q-\overline{Q}|^pd\sigma.
$$
By (\ref{1-5}) and (\ref{1-q}) we have
\begin{equation}\label{g2t}
\partial_tG_2=2\int_0^{2\pi}(Q-\overline{Q})Q_t
d\sigma+\frac12\int_0^{2\pi} (Q-\overline{Q})^3d\sigma,
\end{equation}
and
\begin{equation}\label{qqt}
\int_0^{2\pi}(Q-\overline{Q})Q_t
d\sigma=-\frac{\alpha^2}{12} ||Q_{\sigma\sigma}||_{L^2}^2
+\frac{5\alpha}{6}\int_0^{2\pi}R
Q^2_\sigma d\sigma -2\int_0^{2\pi}Q(Q-\overline{Q})^2d\sigma.
\end{equation}

We have the following well-known interpolation inequality.
\begin{lemma}
$$
\int_0^{2\pi} Q_\sigma^4 d\sigma\le C||Q-\overline{Q}||_{L^2}
||Q_{\sigma\sigma}||_{L^2}^3.
$$
\end{lemma}

Using the above lemma and Young's inequality we obtain that for
any $a>0$:
\begin{align}\label{rqs}
\left|\int_0^{2\pi}RQ^2_\sigma d\sigma\right| \le&||R||_{L^2}
||Q_\sigma||_{L^4}^2 \le C ||R||_{L^2}
||Q-\overline{Q}||_{L^2}^{\frac12}
||Q_{\sigma\sigma}||^{\frac32}\\
\le & C||R||_{L^2}(\frac{1}{4a^4}||Q-\overline{Q}||_{L^2}^2
+\frac{3}{4}a^{\frac43}||Q_{\sigma\sigma}||_{L^2}^2)\nonumber.
\end{align}
Since $||R||_{L^2}^2=2\pi \overline{Q}$ is bounded (see Remark 7
of \cite{NZ1}), it follows from (\ref{g2t}), (\ref{qqt}) and
(\ref{rqs}) that
$$
\partial_tG_2\le
C_1G_2-\frac72\int_0^{2\pi}(Q-\overline{Q})^3d\sigma-C_2
\int_0^{2\pi}Q_{\sigma\sigma}^2d\sigma.
$$
Noticing that
$$
\left|\int_0^{2\pi}(Q-\overline{Q})^3d\sigma\right|\le
||Q-\overline{Q}||_{L^2}^{\frac32}||Q-\overline{Q}||_{L^6}^{\frac32}
\le \frac1{4a^4}||Q-\overline{Q}||_{L^2}^{6}+\frac34a^{\frac43}
||Q-\overline{Q}||_{L^6}^2
$$
and $||Q-\overline{Q}||_{L^6}^2\le C_3
||Q_{\sigma\sigma}||_{L^2}^2$, we obtain that
\begin{equation}\label{dtg2}
\partial_tG_2\le C_4(G_2+G_2^3)-C_5||Q_{\sigma\sigma}||_{L^2}^2.
\end{equation}
\begin{lemma}\label{lem4}
$$
\lim_{t\to\infty}G_2(t)=0,\qquad \int_0^\infty \int_0^{2\pi}
Q^2_{\sigma\sigma}d\sigma dt<\infty.
$$
\end{lemma}
\begin{proof}
From Lemma \ref{lem2} we know that $\int_0^\infty
G_2(t)dt<\infty$. Therefore for any $ \epsilon>0$, there exists
$t_\epsilon>0$, such that $G_2(t_\epsilon) <\epsilon$ and
$\int_{t_\epsilon}^\infty G_2(t)dt<\epsilon$. If
$\epsilon<1/(1+2C_4)$, we must have that $G_2(t)\le 1$ for all
$t>t_\epsilon$. In fact, if not, let $t_*>t_\epsilon$ be the first
time such that $G_2(t_*)=1$. Integrating (\ref{dtg2}) from
$t_\epsilon$ to $t_*$, we obtain that
$$
1-G_2(t_\epsilon)\le C_4(\epsilon+\epsilon),
$$
which implies $\epsilon\ge 1/(1+2C_4)$. Contradiction.
For $\epsilon<1/(1+2C_4)$ and $t>t_\epsilon$, integration (\ref{dtg2}) from
$t_\epsilon$ to $t$, we obtain that
$$
G_2(t)\le G_2(t_\epsilon)+ C_4(\epsilon+\epsilon)\le \epsilon+2C_4\epsilon.
$$
Hence $\lim_{t\to\infty}G_2(t)=0$ and $\int_0^\infty (\int_0^{2\pi}
Q^2_{\sigma\sigma}d\sigma)dt<\infty$.
\end{proof}

 Direct computation yields
\begin{align}\label{qsigmat}
\partial_t||Q_\sigma||_{L^2}^2=-\frac{\alpha^2}6||Q_{\sigma\sigma\sigma}||
_{L^2}^2-&\frac{5\alpha}{3}\int_0^{2\pi}Q_\sigma Q_{\sigma\sigma\sigma}R
d\sigma\\
-&4\overline{Q}||Q_\sigma||_{L^2}^2-\frac{17}2\int_0^{2\pi}Q_\sigma^2
(Q-\overline{Q})d\sigma\nonumber.
\end{align}
For any $t\ge0$, choose $\sigma_0>0$ such that $Q_\sigma(\sigma_0)=0$. Then
$$
|Q_\sigma|=|Q_\sigma-Q_\sigma(\sigma_0)|\le ||Q_{\sigma\sigma}||_{L^1}
\le \sqrt{2\pi}||Q_{\sigma\sigma}||_{L^2},
$$
which implies $||Q_\sigma||_{L^\infty}\le  \sqrt{2\pi}
||Q_{\sigma\sigma}||_{L^2}$. It follows that
\begin{align}\label{qqr}
\left|\int_0^{2\pi}Q_\sigma Q_{\sigma\sigma\sigma}R d\sigma\right|\le&
\frac1C||Q_{\sigma\sigma\sigma}||_{L^2}^2+C\int_0^{2\pi}Q_\sigma^2R^2d
\sigma\\
\le&\frac1C||Q_{\sigma\sigma\sigma}||_{L^2}^2+C ||Q_\sigma||_{L^\infty}^2
||R||_{L^2}^2\nonumber\\
\le& \frac1C||Q_{\sigma\sigma\sigma}||_{L^2}^2+4\pi^2 C
||Q_{\sigma\sigma}||_{L^2}^2\overline{Q}\nonumber,
\end{align}
and
\begin{equation}\label{qq}
\left|\int_0^{2\pi}Q_\sigma^2(Q-\overline{Q})d\sigma \right|
\le||Q_\sigma||_{L^\infty}^2||Q-\overline{Q}||_{L^1}\le 2\pi\sqrt{2\pi}
||Q_{\sigma\sigma}||_{L^2}^2G_2^{\frac12}.
\end{equation}
Substituting (\ref{qqr}) and (\ref{qq}) in (\ref{qsigmat}) and noticing that
$G_2(t)\to 0$, we obtain that
\begin{equation}\label{qst}
\partial_t ||Q_\sigma||_{L^2}^2\le -C_1||Q_{\sigma\sigma\sigma}||_{L^2}^2
+C_2||Q_{\sigma\sigma}||_{L^2}^2\le C_2||Q_{\sigma\sigma}||_{L^2}^2.
\end{equation}
It follows from Lemma \ref{lem4} that for any $\epsilon>0$, there exists
$t_\epsilon>0$, such that
$$
||Q_{\sigma\sigma}||_{L^2}^2(t_\epsilon)<\epsilon\quad\mbox{ and }\quad
\int_{t_\epsilon}^\infty ||Q_{\sigma\sigma}||_{L^2}^2dt<\epsilon.
$$
For any $t>t_\epsilon$, integrating (\ref{qst}) from $t_\epsilon$ to $t$, we
obtain that
$$
||Q_{\sigma}||_{L^2}^2(t)\le ||Q_{\sigma}||_{L^2}^2(t_\epsilon)
+C_2\epsilon\le (2\pi)^2||Q_{\sigma\sigma}||_{L^2}^2(t_\epsilon)+C_2\epsilon
\le C_3\epsilon.
$$
Hence $||Q_{\sigma}||_{L^2}\to 0$ as $t \to +\infty$, which
implies
$$
\lim_{t\to\infty}||Q-\overline{Q}||_{L^\infty}=0.
$$

\section{Exponential convergence of the $4$-$Q$-flow}
We are now ready to derive the exponential convergence for the
metrics under 4-$Q$-curvature flow and thus complete the proof of
Theorem \ref{theorem1}.

 Suppose $g(t)=v^{-\frac43}(\theta,t)g_s$ is a
solution to the flow equation (\ref{normal}) for $\alpha=4$. As in
the proof of Proposition \ref{4b}, for any $t\in[0,\infty)$, we
can choose $\lambda=\lambda(t)>0,\beta=\beta(t)\in[0,2\pi)$ so
that $u(\theta,t):=({\mathcal T}_{\lambda,\beta}v)(\theta)$
satisfies
\begin{equation}\label{oonn}
\int_0^{2\pi}u(t)\cos\theta d\theta=\int_0^{2\pi}u(t)\sin\theta d\theta=0.
\end{equation}
Then $u(\theta,t)$ is uniformly bounded in $H^2(S^1)$ for
$t\in[0,\infty)$. Therefore there exists a sequence $t_n\to
\infty$, such that $u(\theta,t_n) \rightharpoonup
u_\infty(\theta)$ in $H^2(S^1)$. From Sobolev embedding theorem we
have $u(\theta,t_n) \to u_\infty(\theta)$ in $C^{1,a}$ for any
$a\in(0,\frac12)$ and $u\in C^{1, 1/2}$. Since
$$
\int_0^{2\pi} u^{-\frac23}(\theta,t)d\theta =\int_0^{2\pi}
v^{-\frac23}(\theta,t)d\theta=2\pi,
$$
we obtain that $u_\infty(\theta)>0$ and $u_\infty$ satisfies
$$
u_\infty^{\frac53}(\frac{16}9(u_\infty)_{\theta\theta\theta\theta}
+\frac{40}9(u_\infty)_{\theta\theta}+u_\infty)=Q_\infty,
$$
where $Q_\infty=\lim_{t\to\infty}\overline{Q}$. It follows from
(\ref{oonn}) and the classification of solutions of the above ODE
(see the proof of Theorem 4 in \cite{NZ1}) that $u_\infty=1$.
Using the same argument we can prove that any convergent
subsequence of $u(\theta,t)$ converges to $1$. Since $u(\theta,t)$
is uniformly bounded in $H^2(S^1)$, we have $\lim_{t\to
\infty}u(\theta,t)=1$. Hence $Q_\infty=1.$

Define variable $\gamma$ as the inverse of $\theta$ under map
$\omega_{\lambda(t), \beta(t)}$, that is
$\omega_{\lambda,\beta}(\gamma)=\theta$. Noting that $$u(\theta,
t)=v(\omega_{\lambda,\beta}(\theta))\Psi_{\lambda,\beta}(\theta),
 \ \ \ \ \Psi_{\lambda,\beta}(\theta)=\left(
\lambda^2\cos^2\frac{\theta-\beta}2
+\lambda^{-2}\sin^2\frac{\theta-\beta}2\right)^{3/2} $$  and $
\omega_{\lambda,\beta}(\theta)=\beta+\int_\beta^\theta
{\Psi_{\lambda,\beta}}^{-\frac23}d\theta$
 we have
$d\sigma=v(\theta)^{-\frac23}d\theta=u(\gamma)^{-\frac23}d\gamma$.
Since the $4$-curvature and $4$-$Q$ curvature of metric $d \gamma
\otimes d \gamma$ is 1, we obtain that
\begin{align*}
R(\theta(\gamma))&=u(\gamma)(4(u^{\frac13})''(\gamma)+u^{\frac13}(\gamma)),\\
Q(\theta(\gamma))&=u^{\frac53}(\gamma)(\frac{16}9u_{\gamma\gamma\gamma\gamma}
(\gamma)+\frac{40}9u_{\gamma\gamma}(\gamma)+u(\gamma)).
\end{align*}
It follows that $\lim_{t\to\infty}R(\theta)=1$. Therefore
\begin{align}\label{dtg2new}
\partial_tG_2=&-\frac83||Q_{\sigma\sigma}||_{L^2}^2+\frac{20}{3}
\int_0^{2\pi}RQ^2_\sigma d\sigma-4\overline{Q}||Q-\overline{Q}||_{L^2}^2
-\frac72\int_0^{2\pi}(Q-\overline{Q})^3d\sigma\nonumber\\
=&-\frac83||Q_{\sigma\sigma}||_{L^2}^2+(\frac{20}3+o(1))||Q_\sigma||_{L^2}^2
-(4-o(1))||Q-\overline{Q}||_{L^2}^2,
\end{align}
where $o(1) \to 0$ as $t \to +\infty$.

Consider the Fourier series of $Q$:
\begin{equation}\label{fourier}
Q=\overline{Q}+\sum_{n=1}^\infty(a_n\cos(n\sigma)+b_n\sin(n\sigma))
=\tilde{c}+\sum_{n=1}^\infty(\tilde{a}_n\cos(n\gamma)+\tilde{b}_n\sin(n\gamma)).
\end{equation}
Since $u(t)\to 1$ as $t\to\infty$, we obtain that
\begin{equation}\label{ab}
a_n=\tilde{a}_n+o(1)G^{\frac12}_2, \quad
b_n=\tilde{b}_n+o(1)G^{\frac12}_2,\qquad
n=0,1,2,3,\cdots,
\end{equation}
where $a_0=\overline {Q}$ and $\tilde {a}_0=\tilde{c}$. As in
\cite{NZ2} we need estimates on $a_1$ and $b_1$.
\begin{lemma}\label{lemmaon}
For $\phi$ smooth we have
$$
\int_0^{2\pi}(16\phi^{(4)}(\theta)+40\phi''(\theta)+9\phi(\theta))
(\frac23\phi'(\theta)\cos\theta+\phi(\theta)\sin\theta)d\theta=0.
$$
\end{lemma}
\begin{proof}
Integrating by parts we have \begin{equation}\label{temp1}
\int_0^{2\pi}\phi''\phi'''\cos\theta
d\theta=\frac12\int_0^{2\pi}(\phi'')^2 \sin\theta d\theta.
\end{equation}
It follows that
\begin{align*}
&\int_0^{2\pi}(16\phi^{(4)}+40\phi'')(\frac23\phi'\cos\theta+\phi\sin\theta)
d\theta\\
=&\int_0^{2\pi}\phi''(\frac{32}3\phi'''\cos\theta-\frac{16}3\phi''\sin\theta
+48\phi'\cos\theta+24\phi\sin\theta)d\theta\\
=&0+\int_0^{2\pi}(48\phi''\phi'\cos\theta+24\phi''\phi\sin\theta)d\theta\\
=&24\int_0^{2\pi}(\phi')^2\sin\theta d\theta
-24\int_0^{2\pi}(\phi')^2\sin\theta d\theta
-12\int_0^{2\pi}\phi^2\sin\theta d\theta\\
=&-12\int_0^{2\pi}\phi^2\sin\theta d\theta,
\end{align*}
where we used (\ref{temp1}). Also
$$
\int_0^{2\pi}9\phi\cdot
(\frac23\phi'\cos\theta+\phi\sin\theta)d\theta
=12\int_0^{2\pi}\phi^2\sin\theta d\theta.
$$
and the lemma follows.
\end{proof}
The above lemma yields a Kazdan-Warner type identity:
\begin{corollary}\label{cor1}
Given $g=v^{-\frac43}g_s$ with $v$ smooth. Then $4$-$Q$ curvature of $g$
satisfies
\begin{equation}\label{kw}
\int_0^{2\pi}Q_\theta v^{-\frac23}\cos\theta d\theta
=\int_0^{2\pi}Q_\theta v^{-\frac23}\sin\theta d\theta=0.
\end{equation}
\end{corollary}
\begin{proof}
\begin{align*}
\int_0^{2\pi}Q_\theta v^{-\frac23}\cos\theta d\theta
=&-\int_0^{2\pi}Q(-\frac23 v^{-\frac53}v_\theta\cos\theta
-v^{-\frac23}\sin\theta)d\theta\\
=&\frac19\int_0^{2\pi}(16v^{(4)}+40v''+9v)(\frac23v'\cos\theta+v\sin\theta)
d\theta =0.
\end{align*}
Applying the above lemma to $\phi(\theta)=v(\theta+\frac\pi2)$,
we obtain that
$$
\int_0^{2\pi}Q_\theta v^{-\frac23}\sin\theta d\theta=0.
$$
\end{proof}
Using (\ref{kw}) and $\lim_{t\to\infty}u=1$ we have the following
computation
\begin{align*}
\tilde{a}_1=&\frac1\pi\int_0^{2\pi}Q\cos\gamma d\gamma=-\frac1\pi
\int_0^{2\pi}Q_\gamma\sin\gamma d\gamma\\
=&-\frac1\pi\int_0^{2\pi}Q_\gamma\sin\gamma(u^{-\frac23}(\gamma)
-u^{-\frac23}(\gamma)+1) d\gamma\\
=&-\frac1\pi\int_0^{2\pi}Q_\gamma\sin\gamma(-u^{-\frac23}(\gamma)+1)
d\gamma+0\\
=&o(1)\left(\int_0^{2\pi}Q^2_\gamma d\gamma\right)^{\frac12}
=o(1)||Q_\sigma||_{L^2(d\sigma)}.
\end{align*}
Similarly $\tilde{b}_1=o(1)||Q_\sigma||_{L^2}$.
It follws from (\ref{ab}) that $a_1, b_1=o(1)||Q_\sigma||_{L^2}$.

From  Fourier expansion (\ref{fourier}) of $Q$  we obtain that
\begin{align*}
&||Q_{\sigma\sigma}||_{L^2}^2=\pi\sum_{n=1}^\infty n^4(a_n^2+b_n^2),\quad
||Q_{\sigma}||_{L^2}^2=\pi\sum_{n=1}^\infty n^2(a_n^2+b_n^2),\\
&G_2=||Q-\overline{Q}||_{L^2}^2=\pi\sum_{n=1}^\infty (a_n^2+b_n^2).
\end{align*}
It follows from (\ref{dtg2new}) that
$$
\partial_t G_2=\pi\sum_{n=2}^\infty (-\frac83n^4+\frac{20}3n^2-4)(a_n^2+b_n^2)
+o(||Q_\sigma||_{L^2}).
$$
Hence there exists $a>0$, such that $\partial G_2\le -a G_2$, which implies
that
$$
G_2(t)\le Ce^{-at}, \mbox{ for some }C>0.
$$
For any $T>0$ and $\delta\in[0,1]$, integrating (\ref{dtg2})
from $T$ to $T+\delta$ and using the above inequality we have
$$
\int_T^{T+\delta}||Q_{\sigma\sigma}||_{L^2}^2dt \le C_1e^{-aT},
$$
which implies
$$
\int_T^{T+\delta}||Q-\overline{Q}||_{L^\infty}dt\le 2\pi
\int_T^{T+\delta}||Q_{\sigma\sigma}||_{L^2}dt
\le C_2 e^{-\frac a2T}.
$$
Along  $4$-$Q$ flow (\ref{normal}), $v(\theta,t)$ satisfies
$v^{-1}v_t=-\frac34(Q-\overline{Q})$. Integrating from $T$ to
$T+\delta$ we obtain that
$$
|\ln v(\theta,T+\delta)-\ln v(\theta, T)|\le C_3 e^{-\frac a2T},
$$
for any $\theta\in[0,2\pi]$, $T>0$ and $\delta\in[0,1]$. Hence
$\lim_{t\to\infty}v(\theta,t)=v_\infty(\theta)$, with
$||v(t)-v_\infty||_{L^\infty}\le C_4 e^{-\frac a2t}$ and the
$4$-$Q$ curvature of $g_\infty:=v_\infty^{-\frac43}g_s$ is
constant $1$. This completes the proof of Theorem \ref{theorem2}.

\section{Classification of metrics with constant  $1$-$Q$-curvature}
In this section we shall focus on proving Theorem \ref{theorem2}.

Consider the functional $$ {\bf
F}(u)=\frac19\int_0^{2\pi}(u^2_{\theta\theta}-10u^2_\theta+9u^2)d\theta
\left(\int_0^{2\pi}u^{-2/3}(\theta)d\theta\right)^3.
$$
From the proof of Theorem 3 in \cite{NZ1} we know that $$
\inf_{u\in H^2(S^2), \  satisfying (\ref{on}) } F(u)$$ is achieved
by  $v \in H^2(S^2)$, which  satisfies (\ref{on}) and the
Euler-Lagrange equation
\begin{equation}\label{elq}
v_{\theta\theta \theta \theta}+10 v_{\theta\theta } +9v=\tau
v^{-3}
\end{equation}
for a positive constant $\tau$.


Define $V:{\bf R}\to \bf R$ by
$$
V(y)=v(\arctan(2y))(\frac12+2y^2)^{\frac32}.
$$
From conformal invariant properties of $P_g$ (Proposition
\ref{prop0-1}), we obtain that $V(y)$ satisfies
\begin{equation}\label{elqr}
V''''(y)=\tau V(y)^{-\frac53} \ \ \ \ \ \ \mbox{in} \ \ \bf R.
\end{equation}
\begin{lemma}\label{lemmasymm}
Let $w(\theta)=v(\theta)^{\frac13}$. Then $w$ satisfies
\begin{equation}\label{symm}
w^5(\theta)(w''(\theta)+w(\theta))=w^5(\theta+\pi)(w''(\theta+\pi)
+w(\theta+\pi)), \mbox{ for all }\theta \in [0, 2\pi).
\end{equation}
\end{lemma}
\begin{proof}
The first integral of equation (\ref{elqr}) is
$$
V'V'''-\frac12(V'')^2=-\frac32\tau V^{-\frac23}+C.
$$
Since $\lim_{y\to\pm\infty}V(y)=+\infty$, we have
$$
\lim_{y\to\pm\infty }(V'V'''-\frac12(V'')^2 )=C.
$$
Direct computation shows that for $y\to\pm\infty$,
$$V'V'''-\frac12(V'')^2 = 18v^2(\pm\frac\pi2)-4(v'(\pm\frac\pi2))^2
+6v(\pm\frac\pi2)v''(\pm\frac\pi2) + O(y^{-2}).
$$
Since $w^5(w''+w)=v^2-\frac29(v')^2+\frac13vv''$, we know
 that
$$
C=18v^2(\pm\frac\pi2)-4(v'(\pm\frac\pi2))^2+6v(\pm\frac\pi2)v''(\pm\frac\pi2))
=18w^5(\pm\frac\pi2)(w''(\pm\frac\pi2)+w(\pm\frac\pi2)).
$$
By applying an arbitrary shift $\theta\to \theta+\beta$, we obtain
(\ref{symm}).
\end{proof}
Due to intermediate value theorem we may assume without loss of
generality that $v(\frac\pi2)=v(-\frac\pi2)$. The next lemma
indicates that the main difficult in the proof of Theorem
\ref{theorem2} is to match  the derivative of $v$ at north pole
with that at  south pole.
\begin{lemma}\label{ind}
If we also have $v'(\frac\pi2)=v'(-\frac\pi2)$, then
$v^{(k)}(\frac\pi2)=v^{(k)}(-\frac\pi2)$ for all $k\in\bf N$. Moreover
$$
v(\theta)=c\left(\lambda^2\cos^2(\theta-\beta)
+\lambda^{-2}\sin^2(\theta-\beta)\right)^{3/2},
$$
for some $\beta, \lambda,c>0$.
\end{lemma}
\begin{proof}
Let $w=v^{\frac13}$. Since $w'=\frac13v^{-\frac23}v'$, we have $w'(-\frac\pi2)
=w'(\frac\pi2)$. Using Lemma \ref{lemmasymm}, we obtain that
$w''(\frac\pi2)=w''(-\frac\pi2)$. It follows
from $v''=3w^2w''+6w(w')^2$ that $v''(\frac\pi2)=v''(-\frac\pi2)$.
Differentiating (\ref{symm}) and using induction, we obtain that
$v^{(k)}(\frac\pi2)=v^{(k)}(-\frac\pi2)$ for all $k\in\bf N$. It follows that
$g(\theta)=v(\theta/2)$ is smooth on $S^1$. Furthermore $g(\theta)$ satisfies
$$
g''''+\frac52g''+\frac9{16}g=\frac\tau{16}g^{-\frac53}.
$$
Using the same argument as in the last part of the proof of
Theorem 4 in \cite{NZ1}, we obtain that
$$
g(\theta)=c\left(\lambda^2\cos^2\frac{\theta-\beta}2
+\lambda^{-2}\sin^2\frac{\theta-\beta}2\right)^{3/2},
$$
for some $\beta,\lambda,c>0$ and the Lemma follows.
\end{proof}
The proof of Theorem \ref{theorem2} is thus completed if
$v'(\frac\pi2)= v'(-\frac\pi2)$. We are left to consider the case
of  $v'(\frac\pi2)\ne v'(-\frac\pi2)$.

For any $a\in\bf R$, let $V_a(y)=V(y-a)$. Then $V_a(y)$ also
satisfies (\ref{elqr}). It follows that
$\tilde{v}(\theta):=V_a(\frac12\tan\theta)
(2\cos^2\theta)^{\frac32}$ satisfies (\ref{elq}) for
$\theta\in(-\frac\pi2, \frac\pi2)$. Observing that
$$
\tilde{v}(\theta)=V_a(\frac12\tan\theta)
(2\cos^2\theta)^{\frac32}=v(\arctan(\tan\theta-2a))f(\theta),
$$
where $f(\theta):=(\cos^2\theta
+(\sin\theta-2a\cos\theta)^2)^{\frac32}$, we obtain that (recall:
$v(\pi/2)=v(-\pi/2)$)
$$
v_a(\theta):=
\begin{cases}
&v(\arctan(\tan\theta-2a))f(\theta),
\mbox{ when }\theta\in(-\frac\pi2,\frac\pi2)\\
&v(\pi+\arctan(\tan\theta-2a))f(\theta),
\mbox{ when }\theta\in(\frac\pi2,\frac{3\pi}2)\\
&v(\frac\pi2),\mbox{ when }\theta=\pm\frac\pi2
\end{cases}
$$
is smooth on $S^1$ and satisfies (\ref{elq}) and (\ref{on}) (see,
for example, Section 6 in \cite{NZ1} for more details).
Furthermore
$$
v_a''(\pm\frac\pi2)=24a^2v(\pm\frac\pi2)+8av'(\pm\frac\pi2)+v''(\pm\frac\pi2).
$$
Choosing (recall: $v'(\frac\pi2)\ne v'(-\frac\pi2)$)
$$
a=-\frac18\frac{v''(\frac\pi2)-v''(-\frac\pi2)}
{v'(\frac\pi2)-v'(-\frac\pi2)},
$$
we have $v_a''(\frac\pi2)=v_a''(-\frac\pi2)$. Let
$w_a=v_a^{\frac13}$. We have $w_a(\frac\pi2)=w_a(-\frac\pi2)$, and
thus $w_a''(\frac\pi2)=w_a''(-\frac\pi2)$ by Lemma
\ref{lemmasymm}. From $v_a''=3w_a^2w_a''+6w_a(w_a')^2$ we derive
that $v_a'(\frac\pi2)=v_a'(-\frac\pi2)$.   It then follows from
Lemma  \ref{ind} that $v_a^{(k)}(\frac\pi2)=v_a^{(k)}(-\frac\pi2)$
for all $k\in\bf N$ and
$$
v_a(\theta)=c\left(\lambda^2\cos^2(\theta-\beta)
+\lambda^{-2}\sin^2(\theta-\beta)\right)^{3/2},
$$
for some $\lambda, c>0$ and $\beta \in [0, 2\pi)$. Direct
computation shows that $v$ has the form (\ref{form}).  This
completes the proof of Theorem \ref{theorem2}.

\medskip
If $v$ satisfies (\ref{elq}), then $g=v^{-\frac43}g_0$ has
constant $1$-$Q$ curvature. So we  classify  all the constant
$1$-$Q$ curvature metrics on $S^1$ satisfying (\ref{on}).




\section{Exponential convergence of the $1$-$Q$-flow}
Based on the classification result in Theorem \ref{theorem2}, we
shall  prove the exponential convergence of the $1$-$Q$ flow using
a similar argument in Section 4.

Suppose that $g(t)=v^{-\frac43}(\theta,t)g_s$ is a solution to the
flow equation (\ref{normal}) for $\alpha=1$ with initial metric
satisfying (\ref{on}). By Lemma \ref{lemor}, we know that the
metric will  satisfies (\ref{on}) for all $t\ge 0$.  As in the
proof of Proposition \ref{1b}, for any $t\in[0,\infty)$, choose
$\lambda=\lambda(t)>0,\beta=\beta(t)\in[0,2\pi)$ so that
$u(\theta,t):=({\bf T}_{\lambda,\beta}v)(\theta)$ satisfies
$\int_0^{2\pi}u(t)\cos2\theta d\theta=\int_0^{2\pi}u(t)\sin2\theta
d\theta=0$. Here $({\bf
T}_{\lambda,\beta}v)(\theta)=v(\sigma_{\lambda,\beta}(\theta))
\Gamma_{\lambda,\beta}(\theta)$,
$$
\Gamma_{\lambda,\beta}(\theta)=\left(
\lambda^2\cos^2(\theta-\beta)
+\lambda^{-2}\sin^2(\theta-\beta)\right)^{\frac32}
$$
and $\sigma_{\lambda,\beta}(\theta)=\beta+\int_\beta^\theta
{\Gamma_{\lambda,\beta}}^{-\frac23}d\theta$.

As in Section 4, using the classification result (Theorem
\ref{theorem2}) we can prove that $\lim_{t\to\infty}u(\theta,t)=1$
and $Q_\infty:=\lim_{t\to\infty}Q=1$.

Define variable $\gamma$ as the inverse of $\theta$ under map
$\sigma_{\lambda(t), \beta(t)}$, that is
$\sigma_{\lambda,\beta}(\gamma)=\theta$. Then we have
$d\sigma=v(\theta)^{-\frac23}d\theta
=u(\gamma)^{-\frac23}d\gamma$. Since the $1$-curvature and $1$-$Q$
curvature of metric $d\gamma\otimes d\gamma$ is 1, we obtain that
\begin{align*}
R(\theta(\gamma))&=u(\gamma)((u^{\frac13})''(\gamma)+u^{\frac13}(\gamma)),\\
Q(\theta(\gamma))&=u^{\frac53}(\gamma)(\frac{1}9u_{\gamma\gamma\gamma\gamma}
(\gamma)+\frac{10}9u_{\gamma\gamma}(\gamma)+u(\gamma)).
\end{align*}
It follows that $\lim_{t\to\infty}R(\theta)=1$. Therefore
\begin{align}\label{dtg2new1}
\partial_tG_2=&-\frac16||Q_{\sigma\sigma}||_{L^2}^2+\frac53
\int_0^{2\pi}RQ^2_\sigma d\sigma-4\overline{Q}||Q-\overline{Q}||_{L^2}^2
-\frac72\int_0^{2\pi}(Q-\overline{Q})^3d\sigma\nonumber\\
=&-\frac16||Q_{\sigma\sigma}||_{L^2}^2+(\frac53+o(1))||Q_\sigma||_{L^2}^2
-(4-o(1))||Q-\overline{Q}||_{L^2}^2,
\end{align}
where $o(1)\to 0$ as $t \to +\infty$.  Write
\begin{equation}\label{fourier1}
Q=\overline{Q}+\sum_{n=1}^\infty(a_n\cos(n\sigma)+b_n\sin(n\sigma))
=\tilde{c}+\sum_{n=1}^\infty(\tilde{a}_n\cos(n\gamma)+\tilde{b}_n\sin(n\gamma)).
\end{equation}
Since $u(t)\to 1$ as $t\to\infty$, we obtain that
\begin{equation}\label{ab1}
a_n=\tilde{a}_n+o(1)G^{\frac12}_2, \quad
b_n=\tilde{b}_n+o(1)G^{\frac12}_2,\qquad
n=0,1,2,3,\cdots,
\end{equation}
where $a_0=\overline {Q}$ and $\tilde {a}_0=\tilde{c}$. We need
the following lemma to estimate $a_2, \ b_2$.
\begin{lemma}\label{lemmaon1}
For $\phi$ smooth we have
$$
\int_0^{2\pi}(\phi^{(4)}(\theta)+10\phi''(\theta)+9\phi(\theta))
(\frac13\phi'(\theta)\cos(2\theta)+\phi(\theta)\sin(2\theta))d\theta=0.
$$
\end{lemma}
A special case of Lemma \ref{lemmaon1} is the following
Kazdan-Warner type identity.
\begin{corollary}\label{cor2}
Given $g=v^{-\frac43}g_s$ with $v$ smooth. Then $1$-$Q$ curvature of $g$
satisfies
\begin{equation}\label{kw1}
\int_0^{2\pi}Q_\theta v^{-\frac23}\cos(2\theta) d\theta
=\int_0^{2\pi}Q_\theta v^{-\frac23}\sin(2\theta) d\theta=0.
\end{equation}
\end{corollary}
The proofs of Lemma \ref{lemmaon1} and Corollary \ref{cor2} are
very similar to the proofs of Lemma \ref{lemmaon} and Corollary
\ref{cor1}. We shall skip all the details here.

Using (\ref{kw1}) and $\lim_{t\to\infty}u=1$ we have
\begin{align*}
\tilde{a}_2=&\frac1\pi\int_0^{2\pi}Q\cos(2\gamma) d\gamma=-\frac1{2\pi}
\int_0^{2\pi}Q_\gamma\sin(2\gamma) d\gamma\\
=&-\frac1{2\pi}\int_0^{2\pi}Q_\gamma\sin(2\gamma)(u^{-\frac23}(\gamma)
-u^{-\frac23}(\gamma)+1) d\gamma\\
=&-\frac1{2\pi}\int_0^{2\pi}Q_\gamma\sin(2\gamma)(-u^{-\frac23}(\gamma)+1)
d\gamma+0\\
=&o(1)\left(\int_0^{2\pi}Q^2_\gamma d\gamma\right)^{\frac12}
=o(1)||Q_\sigma||_{L^2(d\sigma)}.
\end{align*}
Similarly $\tilde{b}_2=o(1)||Q_\sigma||_{L^2}$. It follws from
(\ref{ab1}) that $a_2, b_2=o(1)||Q_\sigma||_{L^2}$. From the
Fourier expansion (\ref{fourier1}) of $Q$  and (\ref{dtg2new1}) we
obtain that
$$
\partial_t G_2=\pi\sum_{n\in{\bf N}, n\ne 2} (-\frac16n^4+\frac{5}3n^2-4)
(a_n^2+b_n^2)+o(||Q_\sigma||_{L^2}).
$$
Hence there exists $a>0$, such that $\partial G_2\le -a G_2$, which implies
that
$$
G_2(t)\le Ce^{-at}, \mbox{ for some }C>0.
$$
The rest of the proof can be carried out similarly  to the proof
of Theorem \ref{theorem1}. We hereby complete the proof of Theorem
\ref{theorem3}.

\end{document}